\documentclass[11pt]{article}

\usepackage{enumerate}
\usepackage{epsfig}
\usepackage{graphicx}
\usepackage{amssymb}
\usepackage{afterpage}
\usepackage{theorem}
\usepackage{xypic}
\usepackage{psfrag}

\newcommand{\N}{\mathbb{N}}
\newcommand{\Real}{\mathbb{R}}

\newcommand{\C}{\mathbb{C}}

\newcommand{\Z}{\mathbb{Z}}

\newcommand{\D}{\mathbb{D}}

\newcommand{\ov}{\overline}

\newcommand{\OO}{\mbox{$\mathbf{O}$}}
\newcommand{\Otwo}{{\mathbf{O}(2)}}

\newcommand{\proof}{\noindent{\bf Proof:} \quad}
\newcommand{\proofof}[1]{\noindent {\bf Proof of #1} \hspace{0.1in}}
\newcommand{\qed}{\hfill\mbox{\raggedright\rule{0.07in}{0.1in}}\vspace{0.1in}}

\newcommand{\graph}{\mbox{graph}}

\newtheorem{theorem}{Theorem}[section]
\newtheorem{proposition}[theorem]{Proposition}

\newtheorem{lemma}[theorem]{Lemma}

\theorembodyfont{\rmfamily}
\newtheorem{definition}[theorem]{Definition}
\newtheorem{example}[theorem]{Example}
\newtheorem{remarks}[theorem]{Remarks}
\newtheorem{remark}[theorem]{Remark}
\newtheorem{hypo}[theorem]{Assumption}
\newtheorem{question}[theorem]{Question}
\newtheorem{property}[theorem]{}

\newcommand{\Fix}{\mbox{{\rm Fix}}}
\newcommand{\arraystart}{\renewcommand{\arraystretch}{1.5}}
\newcommand{\arrayfinish}{\renewcommand{\arraystretch}{1.2}}

\newcommand{\AND}{\quad\mbox{and}\quad}

\title{On the zero set of $G$-equivariant maps}
\author{P-L. Buono, M. Helmer\\
Faculty of Science\\
University of Ontario Institute of Technology\\
Oshawa, ONT L1H 7K4\\
Canada\\
\\
J.S.W. Lamb\\
Department of Mathematics\\
Imperial College London\\
London SW7 2AZ, UK
}
\begin{document}
\maketitle

\begin{abstract}
Let $G$ be a finite group acting on vector spaces $V$ and $W$ and consider a smooth
$G$-equivariant mapping $f:V\to W$. This paper addresses the question of the zero set near
a zero $x$ of $f$ with isotropy subgroup $G$. It is known from results of Bierstone and Field
on $G$-transversality theory that the zero set in a neighborhood of $x$ is a stratified set.
The purpose of this paper is to partially determine the structure of the stratified set near $x$
using only information from the representations $V$ and $W$. We define an index $s(\Sigma)$
for isotropy subgroups $\Sigma$ of $G$ which is the difference of the dimension of the fixed point subspace
of $\Sigma$ in $V$ and $W$. Our main result states that if $V$ contains a subspace $G$-isomorphic to $W$,
then for every maximal isotropy subgroup $\Sigma$ satisfying $s(\Sigma)>s(G)$, the zero set of $f$ near
$x$ contains a smooth manifold of zeros with isotropy subgroup $\Sigma$ of dimension $s(\Sigma)$.
We also present a systematic method to study the zero sets for group representations $V$ and $W$ which
do not satisfy the conditions of our main theorem. The paper contains many examples and raises several
questions concerning the computation of zero sets of equivariant maps. These results have application
to the bifurcation theory of $G$-reversible equivariant vector fields.
\end{abstract}

\section{Introduction}
The purpose of this paper is to introduce a new perspective and several new results
for the study of $G$-equivariant maps $f:V\to W$ where $V$ and $W$ are possibly
non-isomorphic representations. The main goal of this paper is to investigate
zero sets of such maps in a neighborhood of a zero with full isotropy subgroup $G$.

The widespread appearance of symmetry in differential equation models eventually led to
the establishment of equivariant dynamical system as a subbranch of dynamical systems.
One of the successful uses of equivariant dynamical systems is in the study of bifurcation
problems and equivariant bifurcation theory of vector fields is now a standard tool for the
study of bifurcations in symmetric differential equation models from all areas of science
and engineering.

Fundamental results making up the foundations of the theory of local bifurcations
of vector fields in the presence of symmetry can be found in the works of Michel~\cite{Michel72},
Ruelle~\cite{Ruelle73} and Sattinger~\cite{Sattinger}. The use of singularity theory for studying local
bifurcation problems of equivariant vector fields goes back to Golubitsky and Schaeffer~\cite{GS}
and a comprehensive treatment is found in Golubitsky, Stewart and Schaeffer~\cite{GSS88}. The approach of~\cite{GSS88} has been particularly successful in the study of bifurcation problems arising from mathematical models.
Local bifurcation problems of equivariant vector fields have also been studied using
$G$-transversality theory (a.k.a equivariant general position) and the results can be found
in the work of Field and collaborators, see Field~\cite{Field96,Field07} for details.
For instance, the $G$-transversality approach led to a characterization of the criteria under which
the so-called Maximum Isotropy Subgroup Conjecture (MISC) holds, see for instance~\cite{field-rich1}.

For local bifurcations of $G$-equivariant bifurcation problems $f(x,\lambda)=0$ where $x\in V$,
$\lambda\in\Real^{\ell}$ and $f:V\times\Real^{\ell}\to V$ is $G$-equivariant, the starting point is the
analysis of a zero of the $G$-equivariant mapping $f$. Note that in this case,
the $V$ in the domain and image are isomorphic $G$-representations and the inverse function theorem shows
that equilibrium solutions are generically isolated. However, vector fields which are not only equivariant but also possess properties of time-reversibility, known as $G$-reversible equivariant vector fields, are described
by smooth maps $f:V\times\Real^{\ell}\to V_{\sigma}$ where $V_{\sigma}$ is a $G$-representation of $V$
possibly non-isomorphic to $V$ and determined by the type of time-reversibility. Recent progress on
the steady-state bifurcation theory of $G$-reversible equivariant vector fields~\cite{BLR08} shows
that for large classes of these vector fields, equilibrium solutions are no longer isolated. To determine
the zero set in a neighborhood of equilibrium solutions, one has to study $G$-equivariant maps from
non-isomorphic $G$-spaces $V$ and $W$ where $W$ is a subrepresentation of $V_{\sigma}$.

It is known from $G$-transversality theory that generic zero sets of general $G$-equivariant maps
$f:V\to W$ are Whitney regular stratified sets~\cite{Biers77} and~\cite{Field77}. A stratified set is a locally finite collection of submanifolds and Whitney regularity is a technical condition on the way the submanifolds fit together. Several examples of zero sets are computed in the context of $G$-reversible equivariant systems, see~\cite{Field07} and~\cite{BLR08} where the stratified structure of the zero is partially obtained. The insight gained from these examples shows that partial information about the zero set is encoded in the form of an index which is the difference of the dimensions of fixed point subspaces for isotropy subgroups of $V$ and $W$.

In this paper, we explore this issue and show results which confirm the insight in several cases. We begin by proving in Theorem~\ref{thm:faith1} a nonlinear version of Schur's lemma; that is, we give a sufficient condition on the representations $V$ and $W$ for a $G$-equivariant mapping $f:V\to W$ to be identically zero. Our main theorem is the following: suppose $G$ is a finite group acting on $V$ and $W$ where $V$ contains a subrepresentation
$G$-isomorphic to $W$. Let $f:V\to W$ be a $G$-equivariant map such that $f(0)=0$, then for each maximal isotropy subgroup $\Sigma$ with index greater than the index of $G$, the zero set of $f$ near $0$ contains a submanifold of zeros with isotropy subgroup $\Sigma$ of dimension given by the index. The proof of this result is obtained using a reduction of the problem to isotypic components of $W$ and applying the implicit function theorem. Moreover, for the cases not treated using this result, we present a method suitable for explicit examples and which uses a result of Buchner et al~\cite{BMS83}. This result requires only the computation of the lowest degree equivariants. This is a significant advantage to the alternative method which requires the computation of a minimal set of equivariant generators for smooth $G$-equivariant maps. This is a tedious task, often requiring the use of symbolic algebra packages.

Note that $G$-transversality theory has been formulated in the context of manifolds and that many of the results obtained in this paper can be lifted to smooth mappings between $G$-manifolds using the Slice Theorem~\cite{Field07}.
For instance, $G$-transversality is used in the study of low-dimensional manifolds supporting a group action. In fact, the number which corresponds to the index as defined in this paper appears in Hambleton~\cite{HL92}, but it is not explicitly singled out.

The paper is organized as follows. In the first section, we state and prove Theorem~\ref{thm:faith1} and then present some elementary examples which leads to the statement of our main result, Theorem~\ref{thm:main}. The following section discusses the relevance of these questions in the context of steady-state bifurcations of $G$-reversible equivariant vector fields. Section~\ref{sec:prelim} contains all the prerequisites concerning local zero sets, stratifications and $G$-transversality. Section~\ref{sec:dim-ex} presents known results about the dimension of zero sets with symmetry obtained from stratumwise transversality. Section~\ref{sec:reduc} presents known and new results which enable us to reduce the calculations along isotypic components of the $W$ representation. In section~\ref{sec:inc}, we show Theorem~\ref{thm:ift} and this is the main ingredient in the proof of Theorem~\ref{thm:main} also in this section. Then, we present the computational method based on a result of~\cite{BMS83} to study the cases not covered by Theorem~\ref{thm:ift}. More questions are listed in the final section.

\subsection{Main Theorems}\label{sec:prelim}
Let $V$ and $W$ be finite dimensional vector spaces. Let $G$ be a
compact Lie group and $\rho_{V}:G\rightarrow \mbox{\bf O}(V)$,
$\rho_{W}:G\rightarrow \mbox{\bf O}(W)$ be two representations of
$G$. Recall that if $\rho$ is a representation then
$\ker\rho=\{g\in G\mid \rho(g)=I\}$ and a representation is faithful
if $\ker\rho=\{1\}$.

Let $f:V\rightarrow W$ be a smooth map commuting with the
respective actions of $G$ on $V$ and $W$:
\begin{equation}
f(\rho_{V}(g) x)=\rho_{W}(g)f(x).
\end{equation}
Then $f$ is said to be {\em $G$-equivariant} and we denote this set of functions
by $C_{G}^{\infty}(V,W)$. Let $x\in V$, the set
\[
G_{x}=\{g\in G|gx=x\}
\]
is a subgroup of $G$ called the {\em isotropy subgroup} of $x$. Let
$(G_{x})$ denote the conjugacy class of $G_{x}$, the conjugacy class
of an isotropy subgroup is called the {\em isotropy type}. We write
$\iota(x)$ for the isotropy type of the point $x$. Denote by ${\cal O}(V,G)$
the set of isotropy types for the action of $G$ on $V$. One can define a partial
order on this set by the following rule: let $\tau,\mu\in {\cal O}(V,G)$, then
$\tau>\mu$ if there exists $H\in \tau$ and $K\in \mu$ such that
$H\supsetneq K$. To each isotropy subgroup $\Sigma$ is associated a
{\em fixed point subspace}
\[
\Fix(\Sigma)=\{x\in V| \sigma x=x\;\mbox{for all $\sigma\in
\Sigma$}\}
\]
and for each isotropy type $\tau$, we define the {\em orbit stratum}
$V_{\tau}$ as
\[
V_{\tau}=\{x\in V|\iota(x)=\tau\}.
\]
An important feature of $G$-equivariant maps is that they preserve
fixed point subspaces:
\[
f:\Fix_{V}(\Sigma)\rightarrow \Fix_{W}(\Sigma).
\]
The proof is straightforward as we now show. Let $x\in
\Fix_{V}(\Sigma)$ and $\sigma\in \Sigma$ then
\[
f(x)=f(\rho_{V}(\sigma) x)=\rho_{W}(\sigma)f(x).
\]
We now study the effects of the faithfulness of the representations
$V$ and $W$ on the zero set of $f$. The next result is stated for an
isotypic component $W$. This is not a restriction since the study
of the zero set of $f$ can be decomposed along isotypic components,
more on that in section~\ref{sec:reduc}. Note that this next result can be
interpreted as a nonlinear version of Schur's lemma and treats the case where the
representation $V$ is not faithful and $\ker\rho_V\cap\ker\rho_W=\{1\}$.
\begin{theorem}\label{thm:faith1}
Suppose that assumptions~(\ref{h:trivial}) and~(\ref{h:isotypic}) are satisfied.
Let $(V,\rho_{V})$ and $(W,\rho_{W})$ be $G$-spaces where
$W=U\oplus\cdots\oplus U$ is an isotypic component, and $V$ does
not contain an irreducible representation isomorphic to $U$.
Let $f:V\to W$ be a $G$-equivariant map. Suppose that
$\ker\rho_V\neq \{1\}$ and $\ker\rho_W\cap \ker\rho_V=\{1\}$.
Then, $f$ is identically zero.
\end{theorem}

\proof Suppose that $W$ contains $m$ copies of the irreducible
representation $U$. Let $f(x)=(f_1(x),\ldots,f_m(x))$ where
$f_i:V\to U$ is also $G$-equivariant for $i=1,\ldots,m$ since $U$
is irreducible. Let $\rho_V$ be the representation of $G$ on $V$ and
$\Sigma=\ker\rho_{V}\neq \{1\}$. Then, for all $\sigma\in \Sigma$
\[
f_i(x)=f_i(\sigma x)=\sigma f_i(x)
\]
which means that $f_i(x)\in \Fix_{U}(\Sigma)$ for all $x\in V$. But
since $f_i$ is $G$-equivariant, we have that $\sigma f_i(gx)=f_i(g
x)$ and so $\sigma g f_i(x)=gf_i(x)$ implies
\[
[g^{-1}\sigma g] f_i(x)=f_i(x)
\]
for all $x\in V$. Thus, $$f_i(x)\in \bigcap_{g\in G}
\Fix_U(g^{-1}\Sigma g)$$ for all $x\in V$. We define the subgroup
\[
H:=\langle g^{-1}\sigma g \mid g\in G,\; \sigma\in \Sigma\rangle.
\]
Since $g^{-1}\Sigma g$ is a subgroup of $H$ for all $g\in G$, then
$f_i(x)\in \Fix_{U}(H)$ for all $x\in V$.

For all $k\in G\setminus H$ and $y\in V$, we have $f_i(ky)\in
\Fix_U(H)$ and
\[
f_i(ky)=k f_i(y)\in k \Fix_U(H).
\]
Thus
\[
f_i(x)\in \bigcap_{k\in G} k \Fix_U(H)
\]
since for any $x\in V$ there exists $k\in G$ and $y_k\in V$ such
that $x=k y_k$. But, $\cap_{k\in G} k \Fix_U(H)$ is a $G$-invariant
subspace of the irreducible representation $U$. Therefore,
$\cap_{k\in G} k \Fix_U(H)$ is either $U$ or $\{0\}$. Suppose that
$\cap_{k\in G} k \Fix_U(H)=U$ then $\Fix_U(H)=U$ but this would mean
that $\ker\rho_{W}\supset H\supset \ker\rho_V\neq\{1\}$, implying
$\ker\rho_V \cap\ker\rho_W\neq\{1\}$ which is a contradiction.
Therefore, $\cap_{k\in G} k \Fix_U(H)=\{0\}$, which implies that
$f_i\equiv 0$ for all $i=1,\ldots,m$ and so $f\equiv 0$. \qed

\begin{remark}
Note that a more algebraic proof of this result can be obtained by showing
that the dimension of $G$-equivariant maps of homogeneous degree $d$ is zero for all $d$.
By the trace formula, one can show that
\begin{equation}\label{1}
\dim (\C[V]_{d}\otimes W)^{G}=(\chi_{S^{d}},\chi_{W})
\end{equation}
where $S^{d}$ is the $d^{th}$-symmetric tensor product of $V$,
$\chi$ is the character of the representation and $(\cdot,\cdot)$
is the inner product of characters. We claim that $( \chi_{S^{d}},\chi_{W} )=0$.
Indeed, suppose that $( \chi_{V}^{d},\chi_{U})\neq 0$. Since
there exists a nonidentity element $g\in G$ such that $gv=v$ for all
$v\in V$ then for all integers $d\geq 0$ and for all $w\in
V^{\otimes d}$ we have $gw=w$. Now, $(\chi_{V}^{d},\chi_{U}) \neq 0$ implies
$\chi_{U}(g)=\chi_{U}(1)$ and we have a contradiction since
$g\not\in \ker\rho_U$ by assumption.
\end{remark}
From Theorem~\ref{thm:faith1}, when considering $f^{-1}(0)$, one only needs
to consider the cases $\ker\rho_V\subseteq \ker\rho_W$ and $\ker\rho_W\subseteq
\ker\rho_V$. Suppose that $\ker\rho_W\subsetneq \ker\rho_V$ and
consider $G'=G/\ker\rho_W$. Then, $G'$ acts faithfully on $W$ and
choose $g\in \ker\rho_V\setminus \ker\rho_W$. Then, $\rho_V$
restricted to $g\ker \rho_W$ is the identity matrix. So, $G'$ does
not act faithfully on $V$ and we are in the situation described in
Theorem~\ref{thm:faith1}. Suppose that $\ker\rho_V\subseteq
\ker\rho_W$. Then, $G'=G/\ker\rho_V$ acts faithfully on $V$ and
Theorem~\ref{thm:faith1} does not apply in this case. The next
assumption holds for the remainder of the paper.
\begin{hypo}\label{h:faith}
$\rho_{V}$ is a faithful representation.
\end{hypo}

Zero sets of general $G$-equivariant mappings with given isotropy
subgroup $\Sigma$ typically arise as families where the dimension is
given by the difference in the dimensions of the fixed point subspaces
in $V$ and $W$. The following example illustrates this fact.

\begin{example} Let $V=\Real^{3}$ and $G=\Z_2(R)$ act on
$V=\Real^{3}$ by $R.(x,y,z)=(x,y,-z)$ and $W=\Real^{3}$ as $-R$. This is an
example of a mapping giving rise to a $G$-reversible vector field in $V$. The
general form of the germ at zero is given by:
\[
\dot x=p(x,y,z^2)z,\quad \dot y=q(x,y,z^2)z, \quad \dot
z=r(x,y,z^2)
\]
where $p$, $q$ and $r$ are smooth functions.
The isotropy subgroup $\Z_2$ is such that $\dim\Fix_V(\Z_2)=2$ and
$\dim\Fix_{W}(\Sigma)=1$. Suppose that $(x_0,y_0,0)\in
\Fix_{V}(\Z_2)$ is an equilibrium; that is, $r(x_0,y_0,0)=0$. Let
\[
f(x,y,z)=\left[\begin{array}{c} p(x,y,z^2)z \\
q(x,y,z^2)z \\ r(x,y,z^2)\end{array}\right]
\]
and define
\[
f_{\mathbb{Z}_2}=f|_{\text{Fix}_{V}(\mathbb{Z}_2)}:\Fix_V(\Z_2)\to
\Fix_{W}(\Z_2).
\]
In fact,
\[
f_{\mathbb{Z}_2}(x,y,0)=\left[\begin{array}{c} 0 \\
0 \\ r(x,y,0)\end{array}\right]
\]
and if the nondegeneracy assumption $dr(x_0,y_0,0)\neq 0$ is
satisfied, then the implicit function theorem implies that
$(x_0,y_0,0)$ is part of a one-dimensional family of equilibria with
the same isotropy subgroup. Note that the dimension of the set
of equilibria is the difference between the dimensions of the fixed
point subspaces for $\Z_2$ in $V$ and $W$ and this is a general fact
which is encoded in the "index" definition below. Equilibria with
trivial isotropy subgroup appear as isolated points and are therefore
bounded away from $\Z_2$ symmetric equilibria.
\label{Z2-example}
\end{example}
In Section~\ref{application}, we discuss the results of this paper in the
context of steady-state bifurcations of $G$-reversible equivariant problems.
\begin{definition}
Let $\Sigma$ be a subgroup of $G$. Then the {\em index} of $\Sigma$
is defined as
\[
s(\Sigma;V,W)=\dim\Fix_{V}(\Sigma)-\dim\Fix_{W}(\Sigma).
\]
If there are no ambiguities about the spaces, we can write $s(\Sigma)$.
\end{definition}
This index is defined in Field~\cite{Field07} and also in Buono et al~\cite{BLR08}
in the context of reversible-equivariant vector fields where it is called
the $\sigma$-index where $\sigma$ is the sign-representation of the
group. It is a straightforward application of transversality theory that for finite
groups, given an isotropy subgroup $\Sigma$, generically, zero sets of equilibria form
a smooth $s(\Sigma)$ dimensional manifold, see for instance Field~\cite{Field07}
or section~\ref{sec:dim-ex} below. Example~\ref{Z2-example} shows that
zeros of $f$ with isotropy subgroup $\Z_2$ and $1$ are bounded away from each other.
We now turn to a simple example which illustrates the question of the embedding of zeros
with a given isotropy subgroup $\Sigma$ within zeros with isotropy subgroups
$\Sigma'\subset \Sigma$.
\begin{example}\label{D2ex}
Let $\D_2$ be the group generated by the elements $\kappa$,$\sigma$ which act on
$V=\Real^{2}$ and $W=\Real$ as follows:
\[
\begin{array}{lll}
\kappa.(x,y)=(x,-y), \qquad &  \kappa.u=u,\\
\sigma.(x,y)=(-x,y), \qquad &  \sigma.u=-u.
\end{array}
\]
The isotropy subgroup lattice including the index is given in Figure~\ref{isotlatD2}.
\arraystart
\begin{figure}[ht]
\centerline{%
\xymatrix{ & \D_2 (0) & \\
\Z_{2}(\kappa) (0) \ar[ur] & & \Z_{2}(\sigma) (1) \ar[ul]\\
 & 1 (1) \ar[ul]\ar[ur] & }}
\caption{Isotropy subgroup lattice for the $\D_{2}$ action
on $\Real^2$}
\label{isotlatD2}
\end{figure}
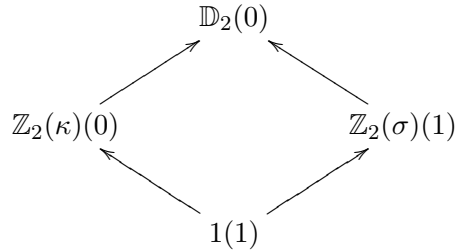
\arrayfinish
Let $f:\Real^2\to \Real$ be a smooth $\D_2$-equivariant mapping.
Since $f:\Fix_{V}(\Sigma)\to \Fix_{W}(\Sigma)$ for all isotropy
subgroups $\Sigma$, zeros of $f$ with isotropy subgroup $\Sigma$
generically appear in $s(\Sigma)$-dimensional families. From this
consideration alone, we see that the zero at the origin with isotropy
subgroup $\D_2$ can be embedded within the family of zeros with
isotropy subgroup $\Z_2(\sigma)$ and $1$, but would be isolated from
zeros with isotropy subgroup $\Z_2(\kappa)$. Zeros with isotropy subgroup
$\Z_2(\kappa)$ could be embedded inside a family of zeros with trivial
isotropy subgroup. We now investigate these options.

The general smooth $\D_2$-equivariant map is
\[
f(x,y)=p(x^{2},y^{2})x.
\]
Therefore $f=0$ if and only if $p(x^2,y^2)=0$ or $x=0$.
$\Fix(\Z_2(\sigma))=\{(x,y)\mid x=0\}$ is a subset of $f^{-1}(0)$
and so the origin is contained inside this one-dimensional family.
A necessary condition for the origin to be embedded within a family of
zeros with trivial isotropy subgroup is that it satisfies the (non-generic)
condition $p(0,0)=0$. Figure~\ref{D2fig}(a) shows this situation. Therefore, the
generic situation in this case is that the origin is embedded within
a family of zeros with isotropy subgroup $\Z_2(\sigma)$, but not with trivial
isotropy subgroup.

Suppose that $p(x_0^2,0)=0$ for some $x_0\neq 0$. The zero at $(x_0,0)$ is included
in a family of zeros with trivial isotropy subgroup if a nondegeneracy condition on
the derivatives of $p(x^2,y^2)$ is satisfied. This generic situation is illustrated in
Figure~\ref{D2fig}(b).
\begin{figure}[htb]
\psfrag{X}{$x$}
\psfrag{Y}{$y$}
\psfrag{1}{$1$}
\psfrag{Z2k}{$\Z_2(\kappa)$}
\psfrag{Z2sigma}{$\Z_2(\sigma)$}
\psfrag{D2}{$\D_2$}
\centerline{%
(a)\includegraphics[height=6.5cm]{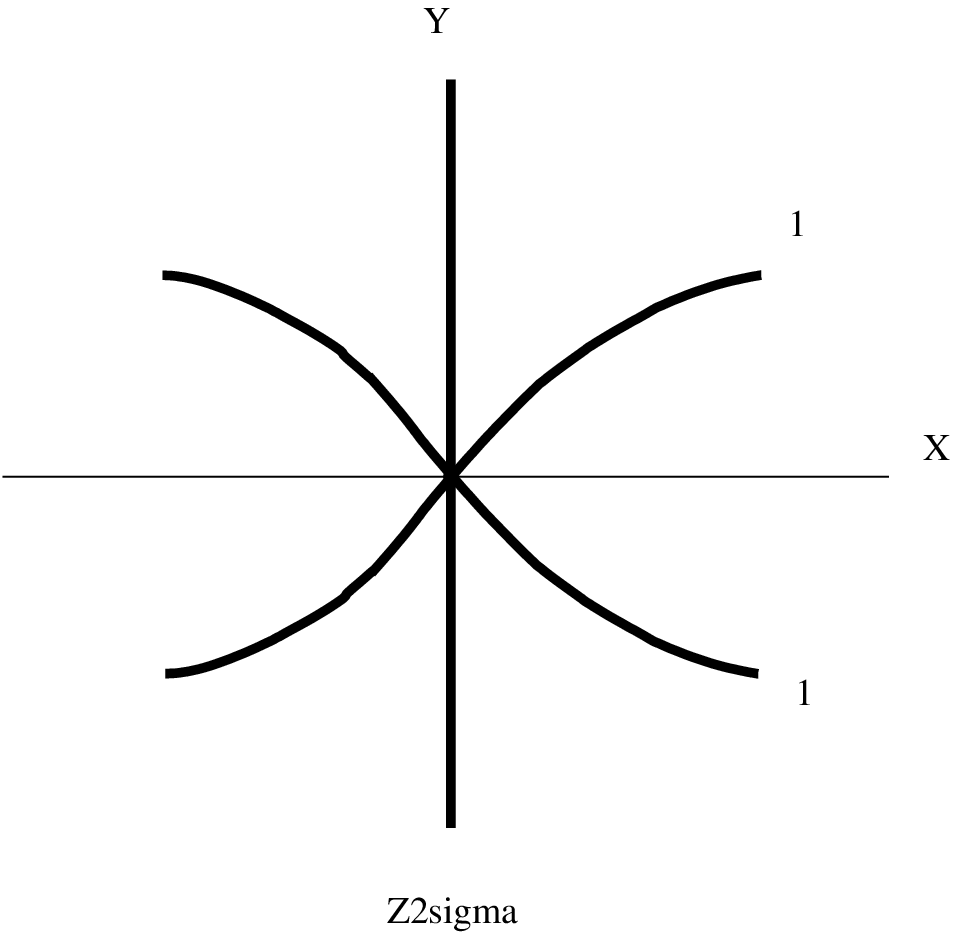}\qquad
(b)\includegraphics[height=6.5cm]{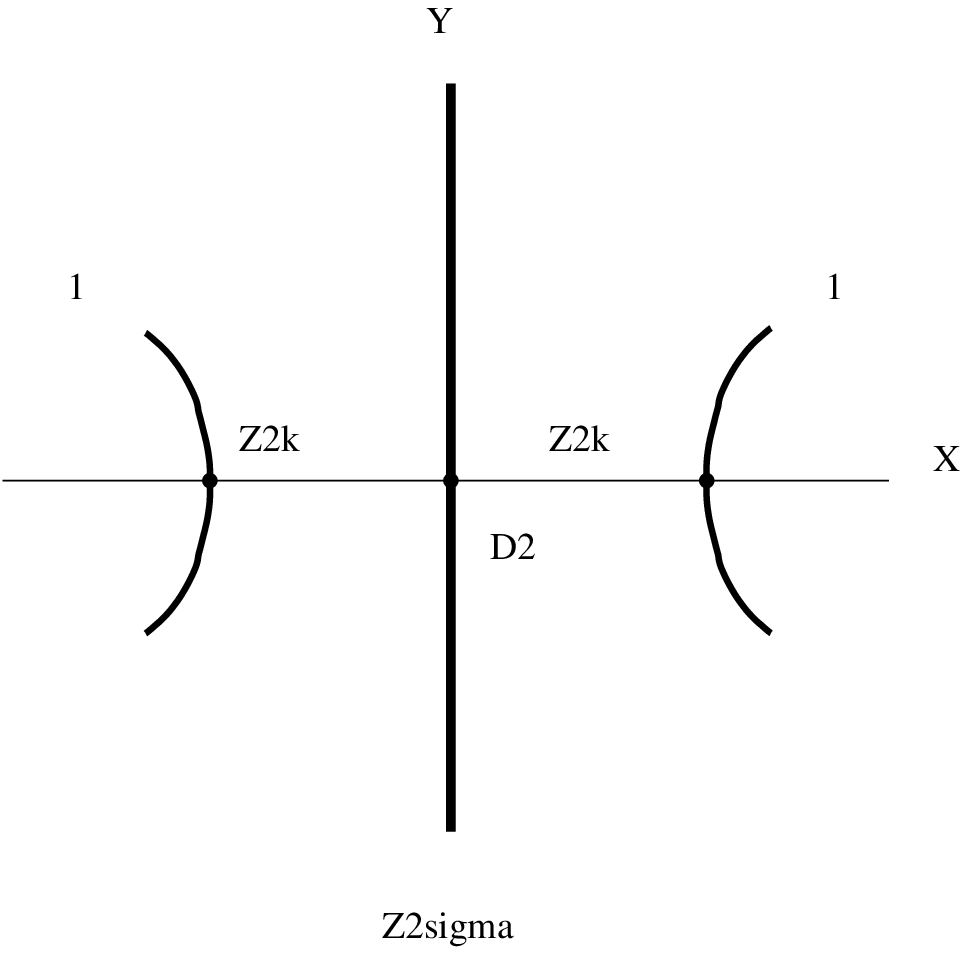}
}
\label{D2fig}
\caption{Zero set of $\D_2$ equivariant map $f$: (a) non-generic case (b) generic case with
$\Z_2(\kappa)$ equilibria}
\end{figure}
\end{example}
Given this example, we can now state our first question.
\begin{question}\label{Q:subgroup}
Let $G$ be a finite group and $\Sigma',\Sigma$ be two isotropy subgroups of $G$. If $\Sigma'$ is a maximal isotropy subgroups of $\Sigma$ and $s(\Sigma')>s(\Sigma)$, is it true that generically the zero set with isotropy subgroup $\Sigma$ is embedded within the zero set with isotropy subgroup $\Sigma'$?
\end{question}
In this paper, we give a partial answer to this question in the case where $\Sigma=G$.
Recall that two representations $V$ and $W$ are $G$-isomorphic if there exists a nonzero
linear $G$-equivariant mapping $A:V\to W$. This is the main result of the paper.
\begin{theorem}\label{thm:main}
Suppose that assumption~\ref{h:faith} is satisfied. Let $f:V\to W$ be a smooth $G$-equivariant map and suppose that $V$
contains a subrepresentation which is $G$-isomorphic to $W$. Suppose that $f$ is $G$-transverse to $0\in W$ at $0\in V$ and let
$\Sigma$ be a maximal isotropy subgroup of $G$ such that $s(\Sigma)>s(G)$. Then, $f^{-1}(0)$ in a neighborhood
of $0\in V$ contains a submanifold of points with isotropy subgroup $\Sigma$ of dimension $s(\Sigma)$
containing $0\in V$.
\end{theorem}
$G$-transversality is defined in section~\ref{sec:strats} and is required in this
result as a nondegeneracy condition.
The proof of Theorem~\ref{thm:main} is a consequence
of the implicit function theorem and it is found in section~\ref{sec:inc}. Now, there are many examples
where the condition on $V$ and $W$ in the above theorem does not hold
and these are studied in section~\ref{sec:case2} and section~\ref{sec:case3}. There, we present an explicit
computational method to determine the structure of the zero set near the origin for representations $V$
and $W$ which do not satisfy the conditions of Theorem~\ref{thm:main}. Moreover, this method only requires
the computation of homogeneous equivariant generators of lowest degree. In section~\ref{sec:questions},
we list some more questions concerning local zero sets of $G$-equivariant mappings.

\section{Application: reversible-equivariant vector fields}\label{application}
The theory developed in this paper has immediate application to the steady-state
bifurcation theory of $G$-reversible equivariant vector fields.
It is a generic feature of reversible-equivariant systems that equilibrium solutions
appear as nontrivial stratified sets.

In~\cite{BLR08}, a systematic study of steady-state bifurcations in smooth
$G$-reversible equivariant vector fields is presented. $G$-reversible equivariant
vector fields are defined as follow. Let $V$ be a vector space and consider the
representations
\begin{equation}\label{rho-sigma}
\begin{array}{rcl}
\rho &: &G\rightarrow \OO(V)\\
\sigma &: &G\rightarrow \Z_{2}=\{+1,-1\}\\
\rho_{\sigma} &: &G\rightarrow \OO(V); \quad
\rho_{\sigma}(g)=\sigma(g)\rho(g).
\end{array}
\end{equation}
A mapping $f$ is $G$-reversible equivariant if for all $g\in G$
\begin{equation}\label{eq:reveq}
f(\rho(g)x,\lambda)=\rho_{\sigma}(g)f(x,\lambda).
\end{equation}
Then, the $\ell$-parameter family of dynamical systems
\begin{equation}\label{eq:f}
\frac{dx}{dt}=f(x,\lambda),
\end{equation}
where $\lambda\in\Real^\ell$ is an $\ell$-dimensional parameter vector
is $G$-reversible equivariant.

Let $L:V\to V$ be a linear $G$-reversible equivariant map commuting
with the representations $\rho$ and $\rho_{\sigma}$ of
$G$ on $V$. If $\rho$ and $\rho_{\sigma}$ are nonisomorphic irreducible
representations, Schur's lemma implies $L\equiv 0$ and
$L$ is forced to have a kernel by the representations. For general
representations $\rho$ and $\rho_{\sigma}$, the map $L$ has a nontrivial
kernel forced by the group representations if the isotypic decompositions
of the representations $\rho$ and $\rho_{\sigma}$ are not isomorphic.
The {\em forced kernel} for a linear reversible equivariant map $L$ is
the isomorphism class (as a group representation) of the lowest dimensional
kernel that a reversible equivariant map can have between these representations.

We decompose $V=W_1\oplus W_2$ where $W_2$ is the forced kernel and $W_1$ is
an orthogonal complement of $W_2$. By construction, $W_2$ and $W_{2,\sigma}$ have
no common irreducible representations and $W_1$ is isomorphic to $W_{1,\sigma}$.
In particular, there exists an invertible $G$-reversible equivariant linear mapping
$T:W_1\to W_{1,\sigma}$. See~\cite{BLR08} for a complete characterization of forced
kernels in terms of irreducible representations of $G$.

Thus, a $G$-reversible equivariant mapping $f:V\to V_{\sigma}$
can be decomposed along $W_1$ and $W_2$ as
\[
f_1:V\to W_1,\quad f_2:V\to W_2.
\]
Suppose that $f(0)=0$ and we want to characterize the zero set in the neighborhood
of $x=0$. The simplest case arises if $d_{W_1}f_1(0)$ is nonsingular, then by the implicit
function theorem, there exists a smooth $G$-reversible equivariant mapping $\phi: W_2\to W_1$
with $\phi(0)=0$ such that in a neighborhood of $w_2=0$,
\[
f_1(\phi(w_2),w_2)\equiv 0.
\]
Therefore, the zero set of $f$ near $x=0$ is characterized by $\tilde{f}_2:W_2\to W_{2,\sigma}$
where $\tilde{f_2}(w_2)=f_2(\phi(w_2),w_2)$ where $W_2$ and $W_{2,\sigma}$ do not share any irreducible
representations. Note that our Theorem~\ref{thm:main} does not hold in this situation and so the
methods presented in sections~\ref{sec:case2} and~\ref{sec:case3} are relevant.
We now return to Example~\ref{Z2-example} to illustrate this general construction.
\begin{example}
The action of $R$ on $V=\Real^{3}$ has two copies of the trivial representation $T$ and one copy of
the $-1$ representation $A$. Let $W_1=T\oplus A=\{(x,y,z)\mid x=0\}$ and $W_2=T=\{(x,y,z)\mid y=z=0\}$.
The mapping $f$ decomposes as $f_1(x,y,z)=(q(x,y,z^2)z,r(x,y,z^2))$ and $f_2(x,y,z)=p(x,y,z^2)z$.
If $f(0)=0$ then $r(0,0,0)=0$ and $d_{W_1}f_1(0,0,0)$ is nonsingular, then $(y,z)=\phi(x):=(\phi_1(x),\phi_2(x))$
where $\phi:W_1\to W_2$ is $R$-equivariant, smooth and $\phi(0)=(0,0)$. Then, $\tilde{f}_2: W_2\to W_{2,\sigma}$
and since $R$ acts trivially on $W_2$ then $\tilde{f}_2(x)=-\tilde{f}_2(x)$. That is, $\tilde{f}_2\equiv 0$ and this
defines a smooth one-dimensional family of zeros of $f$.
\end{example}

At bifurcation points, Buono et al~\cite{BLR08} have shown that $G$-reversible equivariant bifurcation problems
near "organizing centres" can be reduced to $G$-equivariant bifurcation problems with possible
parameter symmetry. Bifurcation problems with parameter symmetry have been studied by~\cite{furter}
and other authors and these bifurcation problems can be analyzed also in the context provided by this paper.

Indeed, consider a smooth vector field $\dot x=f(x,\lambda)$ where $(x,\lambda)\in \Real^{n}\times\Real^{p}$ and suppose that $G$ acts nontrivially on $\Real^{n}$ and the parameters $\Real^{p}$. Then, $G$-equivariance of the vector field means
that $f:V\to W$ where $V=\Real^{n}\times \Real^{p}$ and $W=\Real^{n}$. See Furter {\em et al}~\cite{furter} for more on bifurcation problems with parameter symmetries.

\section{Zero sets and stratifications}\label{sec:strats}
In this section, we introduce concepts and notations in order to study
zero sets of $G$-equivariant maps. A good reference for these results is
the recent book by Field~\cite{Field07}. We discuss briefly the results from stratification
theory which are needed for our purposes and recall the definition of
$G$-transversality.

\paragraph{Local zero sets} The local structure of $f^{-1}(0)$ is investigated
in the following way. The set $C^{\infty}_{G}(V,W)$ is a finitely generated module
over the ring of $G$-invariant maps $C^{\infty}(V)^{G}$. A minimal set of homogeneous
generators (MSG) $\{F_1,\ldots,F_{k}\}$ has the property that all $f\in
C^{\infty}_{G}(V,W)$ can be written as
\[
f=\sum_{i=1}^{k} h_{i}F_{i},
\]
where $h_{i}:V\rightarrow \Real$, $i=1,\ldots,k$ are $G$-invariant
smooth maps. We define the map $F:V\times \Real^{k}\rightarrow
W$ by
\begin{equation}
F(x,t)=\sum_{i=1}^{k} t_{i}F_{i}(x).
\end{equation}
The set ${\cal E}=F^{-1}(0)\subset V\times\Real^{k}$ contains the
information about all possible zero sets for elements $f\in
C^{\infty}_{G}(V,W)$. We refer to ${\cal E}$ as the {\em universal
variety} or also the {\em universal zero set}. The {\em graph map}
of $f$ is a map $\graph_{f}:V\rightarrow V\times\Real^{k}$ defined by
\begin{equation}
\graph_{f}(x)=(x,h_{1}(x),\ldots,h_{k}(x)).
\end{equation}
It is easy to see that
\[
f=F\circ \graph_{f}
\]
and $f^{-1}(0)=\graph_{f}^{-1}({\cal E})$.

\paragraph{Stratifications} We now review some basic facts about
stratifications of sets. Gibson {\em et al}~\cite{GWPL} is a good
reference for the results stated below and more on stratifications.

A {\em stratification} $S$ of a subset $X$ of $\Real^{n}$ is a
locally finite partition of $X$ into smooth and connected
submanifolds of $\Real^{n}$ called {\em strata}. In particular, if
$X$ is semialgebraic, then $S$ is a semialgebraic stratification if
each stratum is semialgebraic. A stratification is said to be {\em
Whitney regular} if it satisfies a Whitney regularity condition. The details
of Whitney regularity are not needed in our work and we refer the interested
reader to~\cite{GWPL}. Of particular interest to us is the {\em frontier condition}
of Whitney regular stratifications: if $A,B$ are strata satisfying $\ov{A}\cap
B\neq \emptyset$ then $B\subset \ov{A}$ or also $B\subset \partial
A$. Another important consequence of Whitney regularity in our context is
the following. If $\ov{A}\cap B\neq \emptyset$, then any map
transverse to $B$ is also transverse to $A$ near $B$. Denote the
union of $i$-dimensional strata of $S$ by $S_{i}$, $i\geq 0$. Then
the following property holds.
\begin{property}\label{inter}
If $S$ is a Whitney regular stratification, then $S_{i}$ does not
meet $\ov{S}_{j}$ unless $i\leq j$.
\end{property}
Now, any semialgebraic set admits a canonical Whitney stratification
having finitely many semialgebraic strata. In particular, let $E$ be
a subset of a smooth manifold $X$. A point $x\in E$ is {\em regular}
(of dimension $d$) if $x$ has a neighborhood $F$ in $X$ such that
$E\cap F$ is a smooth submanifold of (dimension $d$). Every nonempty
semialgebraic set has at least one regular point, in fact the
regular points of maximal dimension are open and dense.

The local structure of zero sets of $G$-equivariant maps has been
determined independently by Bierstone~\cite{Biers77} and Field~\cite{Field77} in their work on
``equivariant general position'' equivalently
``$G$-transversality''.
Let $P$ be a subrepresentation of $W$, then we say that $f$ is $G$-transverse
to $P$ at some point $x$ if $\graph_{f}(x)$ is transverse to the Whitney regular
stratification of ${\cal E}$. In particular, they show the typical properties expected
for transversality theorems; that is, the set
\[
\{f\in C_{G}^{\infty}|\mbox{$f$ is $G$-transverse to $P$}\}
\]
is open and dense in $C_{G}^{\infty}(V,W)$ and an isotopy theorem also holds.

From standard stratification theory one obtains that the structure
of $f^{-1}(P)$ in a neighborhood of a point $x\in f^{-1}(P)$ is an algebraic
set which admits a $G$-invariant Whitney regular stratification. It is sufficient
to study the case $P=\{0\}$ since if $P'$ is a $G$-invariant complement to $P$ in $W$ and
$\pi:W\rightarrow P'$ is the associated projection, then
$f^{-1}(P)=(\pi\circ f)^{-1}(0)$. $G$-transversality of a mapping $f$ to $0\in W$ at $0\in V$ is denoted
by $f\pitchfork_{G} 0\in W$ at $0\in V$.

Hence, the zero set of $f$ is given by the intersection of $\graph_{f}$
with the universal variety ${\cal E}$. The set ${\cal E}$ has a canonical
Whitney regular stratification, thus $\graph_{f}^{-1}({\cal E})$ also has a
Whitney regular stratification, see~\cite{GWPL}.

\section{Dimension of symmetric zero sets}\label{sec:dim-ex}
In this section, we reformulate the problem of the zero set in a
neighborhood of the origin in terms of inclusion of stratum in the
boundary of larger strata. Consider the sets
\[
{\cal E}_{\Sigma}=\{(x,t)\in {\cal E}| x\in V_{\Sigma}\}
\]
over all isotropy types $(\Sigma)$. These are semialgebraic sets and ${\cal E}$ is the disjoint
union of ${\cal E}_{\Sigma}$ over all isotropy types $(\Sigma)$. We denote the union of strata
of regular points (of maximal dimension) in ${\cal E}_{\Sigma}$ by $R_{\Sigma}$. Note that $R_{G}$
is open and dense in ${\cal E}_{G}$. For an isotropy type $(\Sigma)$, a standard notation is $n_{\Sigma}=\dim N_{G}(\Sigma)/\Sigma$.
The next result summarizes important properties of ${\cal E}_{\Sigma}$.
\begin{proposition}[Field~\cite{Field07} Lemma 6.9.2] Let $(\Sigma)$ be an isotropy type for the action of
$G$ and $S$ the Whitney regular stratification of ${\cal E}$. Then,
\begin{enumerate}
\item ${\cal E}_{\Sigma}$ is a semialgebraic submanifold of $V\times\Real^{k}$
of dimension
\[
s(\Sigma;V,W)+\dim (G/\Sigma)-n_{\Sigma}+k,\;\mbox{and}
\]
\item ${\cal E}_{\Sigma}$ inherits a semialgebraic canonical Whitney
stratification $S_{\Sigma}$ from $S$.
\end{enumerate}
\label{prop:dim}
\end{proposition}
For finite groups, the dimension of ${\cal E}_{\Sigma}$ reduces to
$s(\Sigma;V,W)+k$. The intersection property of the sets ${\cal
E}_{\Sigma}$ are given by the next result.
\begin{proposition}[Field~\cite{Field07} Lemma 6.9.1]
Let $(\Sigma'),(\Sigma)$ be isotropy types of $G$. Then
$(\Sigma')>(\Sigma)$ if and only if ${\cal E}_{\Sigma'}\cap
\partial {\cal E}_{\Sigma}\neq \emptyset$.
\end{proposition}
The previous result show that if $(\Sigma')>(\Sigma)$ then some
stratum of ${\cal E}_{\Sigma'}$ is contained in $\partial{\cal
E}_{\Sigma}$. However, to have ${\cal E}_{\Sigma'}\subset \partial{\cal E}_{\Sigma}$
we must show that $R_{\Sigma'}$ is contained in $\partial{\cal E}_{\Sigma}$. The next
result shows a negative criterion for the inclusion of $R_{\Sigma'}$ into $\partial{\cal E}_{\Sigma}$.
\begin{proposition}[Field~\cite{Field07} Lemma 6.9.3]
Suppose that $(\Sigma)<(\Sigma')$. If
$$s(\Sigma';V,W)-n_{\Sigma'}\geq s(\Sigma;V,W)-n_{\Sigma}$$ then
$\dim \left[{\cal E}_{\Sigma'}\cap \ov{\cal E}_{\Sigma}\right]<\dim {\cal
E}_{\Sigma'}$.
\end{proposition}
We now turn to the question of stratumwise transversality and compute the
dimensions of zero sets of $G$-equivariant maps. The next two results are
straightforward consequences of transversality theory, but we state
and prove them explicitly.
\begin{lemma}
Let $f:V\to W$ be $G$--equivariant and suppose that $f(x)=0$, where
$G_{x}=\Sigma$. Then $f\pitchfork_{G} 0$ at $x$ implies that
$s(\Sigma;V,W)-n_{\Sigma}\geq 0$. Moreover, if
$s(\Sigma;V,W)-n_{\Sigma}>0$ then the converse is also true.
\end{lemma}

\proof Since $f\pitchfork_{G} 0$ at $x$, then $\graph_{f}\pitchfork
{\cal E}_{\Sigma}$ at $x$. Now, ${\cal E}_{\Sigma} \subset {\cal
X}=\{G(x,t)|(x,t)\in \Fix_{V}(\Sigma)\times\Real^{k}\}$. So,
$\graph_{f}:\Fix_{V}(\Sigma)\to {\cal X}$ is such that
\[
d(\graph_{f})_{x}(\Fix_{V}(\Sigma))+T_{\graph_{f}(x)}{\cal
E}_{\Sigma} =T_{\graph_{f}(x)}{\cal X}.
\]
Thus,
\begin{equation}\label{eqTr}
\dim d(\graph_{f})_{x}(\Fix_{V}(\Sigma))+\dim T_{\graph_{f}(x)}{\cal
E}_{\Sigma}\geq \dim T_{\graph_{f}(x)}{\cal X}.
\end{equation}
Since $\dim d(\graph_{f})_{x}(\Fix_{V}(\Sigma))=\dim
\Fix_{V}(\Sigma)$, $\dim T_{\graph_{f}(x)}{\cal
E}_{\Sigma}=s(\Sigma;V,W)-n_{\Sigma}+\dim G/\Sigma$ by
Proposition~\ref{prop:dim} and $\dim T_{\graph_{f}(x)}{\cal
X}=\dim \Fix_{V}(\Sigma)+ \dim G/\Sigma+k$, then~(\ref{eqTr})
simplifies to
\[
s(\Sigma;V,W)-n_{\Sigma}\geq 0.
\]
Now, if $s(\Sigma;V,W)-n_{\Sigma}>0$ then the left hand side
of~(\ref{eqTr}) is always greater than the right hand side so
transversality is automatic.\qed

\begin{proposition}\label{redsols}
Suppose that $\graph_{f}$ has nontrivial transverse intersection
with ${\cal E}_{\Sigma}$. Then
\[
\dim \graph_{f}^{-1}({\cal E}_{\Sigma})=s(\Sigma)+\dim (G/\Sigma)-n_{\Sigma}.
\]
In particular, if $G$ is finite then $\dim \graph_{f}^{-1}({\cal E}_{\Sigma})
=s(\Sigma)$.
\end{proposition}

\proof Note that $\dim (d\,\graph_{f})_{x}(V)=\dim V$. Since $\graph_{f}$ has nontrivial
transverse intersection with ${\cal E}_{\Sigma}$, then
\begin{equation}\label{dimeq}
\dim (d\,\graph_{f})_{x}(V)+\dim T_{\graph_{f}(x)}{\cal E}_{\Sigma}
\geq \dim\,V+k.
\end{equation}
which reduces to the dimension of the intersection $s(\Sigma)+\dim
G/\Sigma-n_{\Sigma}\geq 0$.\qed

\begin{remarks}\label{nec-cond}
Note that a necessary condition for the intersection of $\graph_{f}$
with ${\cal E}_{\Sigma}$ to be generically nontrivial is that
\begin{equation}
s(\Sigma)-n_{\Sigma}\geq 0.
\end{equation}
\end{remarks}

\section{Reduction along isotypic components of $W$}\label{sec:reduc}
We begin by a series of results about the structure of the
$G$-spaces $V$ and $W$ which lead to a simplification of the context in
which we study the local zero set near the origin. These results can be
found in Chapter 6 sections 6.6 and 6.7 of Field~\cite{Field07}, but we
include them here for completeness.

Set $p=\dim \Fix_{V}(G)$ and $q=\dim \Fix_{W}(G)$. Let
$V=\Fix_{V}(G)\oplus V'$, $W=\Fix_{W}(G)\oplus W'$ where $V'$ and
$W'$ are the sum of the remaining isotypic components. Let
$\{F_1,\ldots,F_{\ell}\}$ be a MSG for $C_{G}^{\infty}(V,W)$. Let
$\{e_1,\ldots,e_q\}$ be the canonical basis of $\Fix_{W}(G)$. Then
we can set $F_j=e_j$, $j=1,\ldots,q$ and so $F(x,t)=\sum_{i=1}^{q}
t_i e_i+\sum_{i=q+1}^{\ell} t_i F_i(x)$. Therefore $(x,t)\in {\cal
E}$ if and only if $t_j=0$ for $j=1,\ldots,q$ and $\sum_{i=q+1}^{l}
t_i F_i(x)=0$. Write $x=(x_0,x')$ where $x_0\in \Fix_{V}(G)$ and
$x'\in V'$. The coordinates $x_0$ are $G$-invariant so that
$\sum_{i=q+1}^{\ell} t_i F_i(x)=\sum_{i=q+1}^{\ell} t_i F_i(x')$.
Consider $C_{G}^{\infty}(V',W')$ with universal zero set
$\tilde{\cal E}$. We have shown the following correspondence.
\begin{proposition}\label{prop:trivial}
$
(x',t_{q+1},\ldots,t_{\ell})\in \tilde{\cal E}\quad\mbox{if and only
if}\quad (x_0,x',0,\ldots,0,t_{q+1},\ldots,t_{\ell})\in {\cal E}.
$
\end{proposition}
Thus we make the following assumption.
\begin{hypo}\label{h:trivial}
$V$ and $W$ contain no trivial representations.
\end{hypo}
The next two results shows that we can decompose the zero set
problem along isotypic components of the $W$ representation.
\begin{proposition}[Field~\cite{Field07} Lemma 6.6.9]
Let $W=W_1\oplus \cdots \oplus W_{k}$ be the isotypic decomposition
of $W$. Then
\[
\{G_{1}^{1},\ldots,G_{\ell_1}^{1},G_{1}^{2},\ldots,G_{\ell_2}^{2},\ldots,
G_{1}^{k},\ldots,G_{\ell_{k}}^{k}\}
\]
is a MSG for $C_{G}^{\infty}(V,W)$ if and only if for all
$i=1,\ldots,k$, $\{G_{1}^{i},\ldots,G_{\ell_i}^{i}\}$ is a MSG for
$C_{G}^{\infty}(V,W_{i})$.
\end{proposition}
Consider the $G$-equivariant mappings $F^{i}:V\rightarrow W_{i}$ defined
by
\[
F^{i}(x,t^{i})=\sum_{j=1}^{\ell_i} t_{j}^{i} G_{j}^{i}
\]
where $\{G_{1}^{i},\ldots,G_{\ell_i}^{i}\}$ is a MSG for
$C_{G}^{\infty}(V,W_{i})$ and let
\[
{\cal E}_{\Sigma}^{i}=\{(x,t)\in
V_{\Sigma}\times\Real^{\ell_{i}}|F^{i}(x,t)=0\}.
\]
The following proposition shows that inclusion of strata can be obtained
by reducing the inclusion problem to each isotypic component of $W$.
\begin{proposition}\label{prop:ic}
${\cal E}_{G}^{i}\subset \partial {\cal E}_{\Sigma}^{i}$ for all
$i=1,\ldots,k$ if and only if ${\cal E}_{G}\subset \partial {\cal
E}_{\Sigma}$.
\end{proposition}

\proof $\Leftarrow)$ Suppose that ${\cal E}_{G}\subset \partial
{\cal E}_{\Sigma}$. For all $t\in\Real^{l}$, there exists a sequence
$(x^{n},s^{n})\in {\cal E}_{\Sigma}$ such that
$(x^{n},s^{n})\rightarrow (0,t)$. Write $t=(t_1,\ldots,t_k)$ where
$t_{i}\in\Real^{\ell_{i}}$ and for all $n\in\N$, $s^{n}=(s_1^{n},\ldots,s_{k}^{n})$
where $s_{i}^{n}\in \Real^{\ell_i}$. By definition of $F^{i}$, $(x^{n},s_{i}^{n})\in {\cal
E}_{\Sigma}^{i}$ for all $n$. Hence for any $t_{i}\in
\Real^{l_{i}}$, there exists $(x^{n},s_{i}^{n})\rightarrow (0,t_i)$;
that is, ${\cal E}_{G}^{i}\subset \partial {\cal E}_{\Sigma}^{i}$.

$\Rightarrow)$ Suppose that ${\cal E}_{G}^{i}\subset \partial {\cal
E}_{\Sigma}^{i}$ for all $i$. Choose any $t\in \Real^{l}$,
$t=(t_1,\ldots,t_k)$ where $t_{i}\in\Real^{l_{i}}$ and let
$s^{n}=(s_1^{n},\ldots,s_k^{n})$ where $s_{i}^{n}\in \Real^{\ell_i}$ and
$(x^{n},s_{i}^{n})\in {\cal E}_{\Sigma}^{i}$. For all $n$,
$F(x^{n},s^{n})=F^1(x^{n},s_1^{n})+ \ldots+F^k(x^{n},s_k^{n})=0$.
Thus $(x^{n},s^{n})\in {\cal E}_{\Sigma}$ and $(x^n,s^n)\rightarrow
(0,t)$. \qed

For the remainder of the paper we make the next assumption.
\begin{hypo}\label{h:isotypic}
$W$ is the direct sum of $r$ irreducible representations isomorphic
to $U$.
\end{hypo}

\noindent In the next section, we show the main results concerning the inclusion of
${\cal E}_{G}$ inside $\partial {\cal E}_{\Sigma}$.

\section{Inclusions of ${\cal E}_{G}$}\label{sec:inc}
As mentioned above, we assume for the remainder of the paper that
Assumptions~\ref{h:trivial} and~\ref{h:isotypic} are
in force. The MSG for $C_{G}^{\infty}(V,W)$ is taken to be $\{F_1,\ldots,F_{\ell}\}$.
The following remark summarizes the approach taken to obtain
the proofs of the theorems.

\begin{remark}
Recall that $R_{G}$ is the subset of regular points (of
maximal dimension) of ${\cal E}_{G}$ and that $R_{G}$ is open and dense
in ${\cal E}_{G}$. Therefore, the inclusion ${\cal E}_{G}\subset \partial {\cal
E}_{\Sigma}$ is equivalent to $R_{G}\subset \partial {\cal
E}_{\Sigma}$. Now, since ${\cal E}$ has a
Whitney regular stratification, by the frontier condition, to prove the inclusion $R_{G}\subset
\partial {\cal E}_{\Sigma}$ it is sufficient to show that there exists one element $(0,t)\in R_{G}$
and a sequence $\{(x_n,t^{n})\}\subset {\cal E}_{\Sigma}$ such that $(x_n,t^{n})$ converges to $(0,t)$.
\end{remark}
It is straightforward that if $\Sigma$ is a maximal isotropy subgroup of $G$ with
$s(\Sigma)>0$ and $\Fix_{W}(\Sigma)=\{0\}$ then ${\cal E}_{G}\subset
\partial {\cal E}_{\Sigma}$. We suppose hereafter that $\Fix_{W}(\Sigma)\neq \{0\}$.
Suppose that $W$ is the direct sum of $r$ irreducible representations $U$ and let $\delta:=(\chi_{V},\chi_{U})$
where $\chi_{V}$ and $\chi_{U}$ are the characters of $V$ and $U$. Three cases have to be considered:
\begin{enumerate}
\item $\delta\geq r$
\item $\delta=0$, and
\item $0<\delta<r$.
\end{enumerate}
The first case can be solved completely and is the topic of the next section.
The second and third cases are not as straightforward and we do not present a
general theorem. However, we show a systematic method to study particular examples
in those cases.

\subsection{First Case}
This case is the only one which is amenable completely to the implicit function theorem.
Therefore, we obtain a general result on the inclusion of ${\cal E}_G$ into the boundary
of neighbouring sets ${\cal E}_{\Sigma}$.
\begin{theorem}\label{thm:ift}
Let $\Sigma$ be a maximal isotropy subgroup of $G$ and suppose
that $s(\Sigma)>0$. If $\delta\geq r$, then ${\cal E}_{G}\subset
\partial {\cal E}_{\Sigma}$
\end{theorem}

\proof By Schur's lemma, if $L:V\rightarrow W$ is a linear map, then
\[
L=[A,0]
\]
where $A:\delta U\rightarrow W$ is
\[
A=\left[\begin{array}{ccc}
a_{11}I & \cdots & a_{1\delta}I\\
\vdots & \ddots & \vdots\\
a_{r1}I & \cdots & a_{r\delta}I
\end{array}\right]
\]
for $a_{ij}\in\mathbf{k}$, $\mathbf{k}\in \Real,\C$.

By the above, there are $r\delta$ linear generators in the MSG for
$C_{G}^{\infty}(V,W)$. Let
\[
L_0=d_{x}F(0,t)=\sum_{i=1}^{r\delta} t_i L_i.
\]
The matrix $(a_{ij})$ is generically of rank $r$. Hence there exists
$\tau\in R_{G}$ such that for $(\tau_1,\ldots,\tau_{r\delta})$,
$L_0$ is of rank: $r\times\dim U$.

Let $V=rU\oplus V'$ where $V'$ is a $G$-invariant complement of a subspace which is the
direct sum of $r$ copies of $U$. By the implicit function theorem, there exists open neighborhoods
$N_1$ of $0\in V'$ and $N_2$ of $0\in rU$ and a smooth $G$-equivariant mapping $\phi:N_1\to N_2$ such that
\[
F(\phi(v'),v',\tau)\equiv 0.
\]
Let $\Sigma$ be a maximal isotropy subgroup of $G$ with $s(\Sigma)>0$.
Then, $\dim\Fix_{V'}(\Sigma)=s(\Sigma)$ and $(v',\phi(v'))$ for $v'\in N_1$
defines a smooth $s(\Sigma)$-dimensional manifold through the origin in $\Fix_{V}(\Sigma)$.
Therefore, ${\cal E}_{G}\subset \partial {\cal E}_{\Sigma}$.\qed

The next example illustrates the use of Theorem~\ref{thm:ift}. It is not standard
in the equivariant bifurcation theory literature, so we present it in details. The group
under study contains the smallest group in the family $F_{p,q}$ with $p$ prime and $q\mid p-1$
for which there are non-isomorphic representations with maximal isotropy subgroups having
dimension two; in fact, there are three such non-isomorphic representations. See~\cite{JL}
for a description of the family $F_{p,q}$.

\begin{example}\label{ex:F-group}
Consider the group $\Gamma=F_{13,4}\times \mathbb{Z}_{2} $,
where
\[
F_{13,4}=\left\langle a, b:
a^{13}=b^{4}=1,b^{-1}ab=a^{5}\right\rangle.
\]
$F_{13,4}$ has order $52$, so $\Gamma$ has order 104. Let $\omega = e^{\frac{2\pi i}{13}}$, and
\[
\alpha= \omega + \omega^{5} +\omega ^{8} +\omega^{12}, \quad \beta =
\omega^{2} + \omega^{3} +\omega ^{10} +\omega^{11}, \quad \gamma =
\omega^{4} + \omega^{6} +\omega ^{7} +\omega^{9}.
\]
\arraystart
\begin{table}
\begin{center}
\begin{tabular}{ l | c c c c c c c c c c c}
\hline \hline
$g_{i}$ & 1 &-1 & $a$ & $a^{2}$ & $a^{4}$ & $ \pm b$ & $\pm b^{2}$ & $ \pm b^{3}$& $-a$ & $-a^{2}$ & $-a^{4}$\\
$|C_{G}(g_{i})|$ & $52$ & $52$ &$13$ & $13$ & $13$ & $4$ & $4$ & $4$& $13$ & $13$ & $13$\\
\hline
   $\phi_{1}$& 4 &-4&$\alpha$ &$\beta$ &$\gamma$& 0 & 0 &0&$-\alpha$ &$-\beta$ &$-\gamma$\\
   $\phi_{2}$& 4 &-4&$\beta$ &$\gamma$ &$\alpha$& 0 & 0 &0&$-\beta$ &$-\gamma$ &$-\alpha$\\
   $\phi_{3}$& 4 &-4&$\gamma$ &$\alpha$ &$\beta$& 0 & 0 &0&$-\gamma$ &$-\alpha$ &$-\beta$\\
\hline
\end{tabular}
\end{center}
\caption{Some characters of irreducible representations for $F_{13,4}$.}
\label{tab:char}
\end{table}
\arrayfinish
The character table in Figure~\ref{tab:char} shows the characters of the
three irreducible representations of degree $4$, the full character table
of $F_{13,4}$ is in~\cite{JL}. These three irreducible representations are
denoted $V_{1}$, $V_{2}$ and $V_{3}$. For $\mu,\nu\in V_{i}$,
$i\in\left\{1,2,3\right\}$,
\[
b.(\mu,\nu)=(\overline{\nu},\mu)
\]
and for the $a$-actions we have for $(z_{1}, z_{2})\in
V_{1}$, $(z_{3}, z_{4})\in V_{2}$ and $(z_{5}, z_{6})\in V_{3}$
\[
a . (z_{1}, z_{2})=(\omega z_{1}, \omega^{5} z_{2})\quad
a . (z_{3}, z_{4})=(\omega^{10} z_{3},\omega^{11} z_{4})\AND
a . (z_{5},z_{6})=(\omega^{7} z_{5}, \omega^{9} z_{6}).
\]
We suppose that $\Z_2$ acts as $-I$ on those representations.
The lattice of isotropy subgroups is given in
Figure~\ref{latticeF134}.
\begin{figure}[ht]
\centerline{%
\xymatrix{ & \Gamma  & \\
\mathbb{Z}_{4} ( b )\ar[ur] &
\mathbb{Z}_{4}( -b)\ar[u]&\mathbb{Z}_{2}(-b^{2} )\ar[ul]\\
\mathbb{Z}_{2}( b^{2})\ar[ur]\ar[u]\\
  & \ar[ul] 1 \ar[uu]\ar[uur]& }}
\caption{Isotropy subgroup lattice for $\Gamma$ action on
$V_{1},V_{3},V_{3}$} \label{latticeF134}
\end{figure}
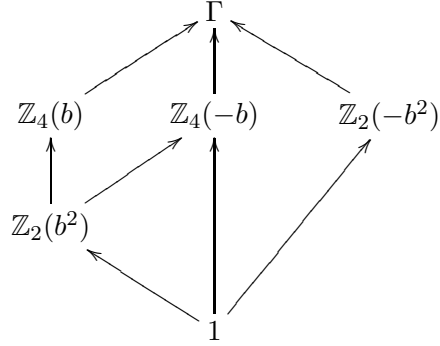
Let $\mu= x+ i y$ and $\nu=u+iv$, with $x,y,u,v \in \Real$. For $j
\in \left\{1,2,3\right\}$ we have
\[
\Fix_{V_{j}}(\Z_4(b))=\left\{(x,0,x,0)\right\},
\quad
\Fix_{V_{j}}(\Z_4(-b))=\left\{(x,0,-x,0)\right\},
\]
and both fixed point subspaces have dimension one.
Moreover,
\[
\Fix_{V_{j}}(\Z_2(b^{2}))=\left\{(x,0,u,0)\right\},\quad
\Fix_{V_{j}}(\Z_2(-b^{2}))=\left\{(0,y,0,v)\right\}
\]
and so both fixed point subspaces have dimension two. Thus, $\Z_2(-b^{2})$ is a maximal
isotropy subgroup with fixed point subspace of dimension $2$.

Considering mappings of the type $f:V_{i} \times V_{j} \rightarrow
V_{k}$, with $i, j, k \in \left\{1, 2, 3\right\}$ then the indices
of the isotropy subgroups are:
\[
s(\Z_4(b))=1, \quad s(\Z_4(-b))=1 ,\quad s(\Z_2(b^{2}))=2,
\quad s(\Z_2(-b^{2}))=2.
\]
Letting $f:V_{1}\times V_{2} \rightarrow V_{2}$ the lowest degree
equivariant is the linear equivariant, since $V_{2}$ is present in
both the domain and range of $f$, this equivariant is a
scalar multiple of the identity endomorphism from $V_{2}$ to $V_{2}$
by Schur's Lemma. Looking at Theorem~\ref{thm:ift}, we have $\delta =1$ as
the domain and image share one irreducible representation; also
there is only one isotypic component in the range therefore the
assumptions of Theorem~\ref{thm:ift} are satisfied and we can conclude that
${\cal E}_{G}\subset\partial {\cal E}_{\Sigma}$ for the maximal
isotropy subgroups $\Z_4(b)$, $\Z_4(-b)$ and $\Z_2(-b^{2})$.
\end{example}

\noindent Now, using Theorem~\ref{thm:ift}, we can finally prove our main result.

\vspace{1\baselineskip}

\proofof{Theorem~\ref{thm:main}} From Proposition~\ref{prop:trivial}, we neglect
the trivial representations of $V$ and $W$ for studying the zero set near $0\in V$.
Decompose the study of $f^{-1}(0)$ along isotypic components
of $W$. That is, let $f=(f_1,\ldots,f_k)$ where $f_{i}:V\to W_i$ and consider $f_{i}^{-1}(0)$
for all $i=1\ldots,k$. The representations $V$ and $W_i$ satisfy the hypothesis
of Theorem~\ref{thm:ift}, therefore for all maximal isotropy types $(\Sigma)$ of $G$
with $s(\Sigma)>s(G)$, ${\cal E}_{G}^{i}\subset \partial {\cal E}_{\Sigma}^{i}$ for
all $i=1,\ldots,k$. By Proposition~\ref{prop:ic}, ${\cal E}_{G}\subset \partial{\cal E}_{\Sigma}$
for all the above mentioned isotropy types. Since $f$ is $G$-transverse to
$0\in W$ at $0\in V$, $\graph_{f}$ intersects transversally a stratum $S_G$ of ${\cal E}_G$.
Then it also intersects nontrivially and transversally a stratum of ${\cal E}_{\Sigma}$ for all
maximal isotropy types $(\Sigma)$ of $G$ with $s(\Sigma)>s(G)$. Hence,
$f^{-1}(0)=\graph_{f}^{-1}({\cal E})$ near $0\in V$ contains branches of zeros for all
maximal isotropy subgroups $\Sigma$ of $G$ with $s(\Sigma)>s(G)$ and these have dimension
$s(\Sigma)$ by Proposition~\ref{redsols}.\qed

\subsection{Second case}\label{sec:case2}
We now turn to the case where $\delta=0$. That is, $V$ and $W$ do not share irreducible
representations and so the smallest degree of nonzero $G$-equivariant homogeneous generators
is $d\geq 2$. Hence, we cannot use the implicit function theorem as
in the previous case. Instead, we show below how a theorem of Buchner {\em et al}~\cite{BMS83}
can be used systematically to study the inclusion of ${\cal E}_{G}$ into $\partial {\cal E}_{\Sigma}$
when $\Sigma$ is a maximal isotropy subgroup with positive index.

A homogeneous polynomial $Q(x)$ is said to be {\em regular on its zero set} if
the Jacobian matrix $dQ(x)$ is surjective for all $x\in Q^{-1}(0)\setminus \{0\}$.
Here is the required result.
\begin{theorem}[Buchner et al~\cite{BMS83}]\label{BMS}
Let $k\geq 2$ be an integer. Suppose $g:\Real^{n}\rightarrow \Real^{m}$ is
smooth, $\Gamma$-equivariant and $g(0)=0$, $dg(0)=0,\ldots,d^{k-1}g(0)=0$.
Let $Q$ be the $k$-form associated to $d^{k}g(0)$ and assume that $Q$ is
regular on its zero set. Then there are $\Gamma$-invariant neighborhoods
$U_1$, $U_2$ containing $0\in \Real^{n}$ and a smooth $\Gamma$-equivariant
diffeomorphism $\phi:U_1\rightarrow U_2$ such that
\[
\phi(Q^{-1}(0)\cap U_1)=g^{-1}(0)\cap U_2,\quad \phi(0)=0
\quad\mbox{and}\quad d\phi(0)=\mbox{identity}.
\]
\end{theorem}
The important assumption to satisfy in order to use Theorem~\ref{BMS}
is the regularity of $Q$ on its zero set. In fact we need regularity
on the zero set only for certain isotropy subgroups.
\begin{lemma}\label{lemma:Q2}
If $Q^{-1}(0)\neq \{0\}$ for some maximal isotropy subgroup
$\Sigma$ with $s(\Sigma)\leq 0$, then $Q$ is not regular on its
zero set.
\end{lemma}
\proof This is obvious for $s(\Sigma)<0$ since $\dim\Fix_{V}(\Sigma)$ is the maximal rank of
$dQ$ restricted to $\Fix_{V}(\Sigma)$ and $$\dim\Fix_{V}(\Sigma)<
\dim\Fix_{W}(\Sigma).$$ If $s(\Sigma)=0$ and $Q(x,\tau)=0$ for some
$x\neq 0$, then $dQ|\Fix_{V}(\Sigma)$ has a nonzero kernel since the
line through $x$ consists of zeroes of $Q$. \qed

We now illustrate the use of Theorem~\ref{BMS} with a few examples.
\begin{example}\label{ex:D6}
Consider the group $G=\D_6$ acting on $V=\C^{2}$ and $W=\C$ as follows.
Let $(z_1,z_2)\in V$ and $w\in W$, then
\[
\kappa.(z_1,z_2)=(\ov{z}_1,\ov{z}_2),\quad \sigma.(z_1,z_2)=(e^{i\pi/3}z_1,
e^{i\pi/3}z_2)
\]
and
\[
\kappa.w=\ov{w},\quad \sigma.w=e^{2i\pi/3}w.
\]
The lattice of isotropy subgroups with indices is given in Figure~\ref{fig:D6}.
We see that both $\Z_2(\kappa)$ and $\Z_2(\kappa\sigma)$ have positive index
so that the zero set near $0$ may include points with those isotropy types.
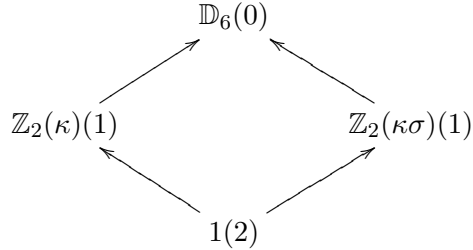
\begin{figure}[ht]
\centerline{%
\xymatrix{ & \D_6 (0)& \\
\Z_{2}(\kappa) (1)\ar[ur] & & \Z_{2}(\kappa\sigma) (1)\ar[ul]\\
 & 1 (2)\ar[ul]\ar[ur] & } }
 \caption{Lattice of isotropy types with index}
 \label{fig:D6}
\end{figure}
The universal zero map is given by
\[
F(z_1,z_2,t)=t_1 z_1^{2}+t_2 z_1 z_2+t_3 z_2^{2}+t_4\ov{z}_1^{4}+t_5\ov{z}_1^{3}
\ov{z}_2+t_6\ov{z}_1^{2}\ov{z}_2^{2}+t_7 \ov{z}_1 \ov{z}_2^{3}
+t_8 \ov{z}_2^{4}
\]
where $t\in \Real^{8}$. The mapping $Q$ of lowest order equivariant terms is given by
\[
Q(z_1,z_2,t_1,t_2,t_3)=t_1 z_1^{2}+t_2 z_1 z_2+t_3 z_2^{2}
\]
and we look at the restriction of $Q$ on the fixed point subspaces of $\Z_2(\kappa)$ and $\Z_2(\kappa\sigma)$.
We have
\[
\Fix_{V}(\Z_2(\kappa))=\{(z_1,z_2)\mid z_1=u\in \Real, \;z_2=v\in \Real\}
\]
so that $Q^{\Sigma}=Q|_{\Fix_{V}(\Z_2(\kappa))}=t_1 u^2+t_2 uv +t_3 v^2$. Thus,
\[
dQ^{\Sigma}(u,v)=[2 t_1 u+t_2 v,t_2 u+2 t_3 v].
\]
We now show that $Q$ is regular on its zero set for an open and dense set of values of $t$.
The system $dQ^{\Sigma}(u,v)=0$ has a nonzero solution if and only if $\Delta=4t_1 t_3-t_2^2=0$. This means that
$Q$ is regular on its zero set for all values of $(t_1,t_2,t_3)$ in the set $\Real^{3}\setminus \Delta$; that is, for an open and dense set in $\Real^{3}$. Now, we can choose $(t_1^{*},t_2^{*},t_3^{*})\in [\Real^{3}\setminus \Delta]$ such that $t^{*}\in R_{G}$ and there exists $(u,v)\neq (0,0)$ such that $Q^{\Sigma}(u,v,t_1^{*},t_2^{*},t_3^{*})=0$. By homogeneity of $Q^{\Sigma}$, we know that $Q^{\Sigma}(\lambda u,\lambda v,t_1^{*},t_2^{*},t_3^{*})=0$ for all $\lambda\in \Real$; that is, there is a line of zeros of $Q^{\Sigma}$ through the origin. Now, consider $F(u,v,t^{*})$. By Theorem~\ref{BMS}, $Q^{-1}(0)|_{\Fix_{V}(\Sigma)}$ is homeomorphic to $F^{-1}(0)|_{\Fix_{V}(\Sigma)}$ near $0$. Hence, there exists a sequence ${(x_n,t^{*})}\subset {\cal E}_{\Sigma}$ such that $(x_n,t^{*})\rightarrow (0,t^{*})$. That is, ${\cal E}_{G}\subset \partial {\cal E}_{\Z_2(\kappa)}$. A similar argument holds for the isotropy subgroup $\Z_2(\kappa\sigma)$. Thus, if a map in $C_{G}^{\infty}(\C^{2},\C)$ is $G$-transverse to $0\in \C$ at $0\in \C^2$, then
the zero set in a neighborhood of the origin contains smooth $1$-dimensional submanifolds of points with isotropy types
$\Z_2(\kappa)$ and $\Z_2(\kappa\sigma)$.
\end{example}

\begin{example}
We return to the case of Example~\ref{ex:F-group} and we now consider the case of $C_{G}^{\infty}(V_{1}\times V_{2},V_{3})$. The polynomial equivariant of smallest degree is $3$ since $\delta=(\chi_{V_1\times V_2},\chi_{V_3})=0$ and there can be no quadratic equivariants because of the nontrivial action of $-I$. In order to apply Theorem~\ref{BMS}, we need to compute the cubic equivariants of the universal zero map $F$. The space of cubic equivariants for $F$ is nine-dimensional and the mapping $Q$ of cubic equivariants is
\[
\begin{array}{rcl}
Q(z,t)&=&t_{1}\left[\begin{array}{cc}z_{4}^{3}\\ \overline{z}_{3}^{3} \end{array}\right]+t_{2}\left[\begin{array}{cc}\overline{z}_{3}\overline{z}_{4}^{2}\\
 {z}_{3}^{2}\overline{z}_{4}\end{array} \right]+t_{3}\left[\begin{array}{cc}
 \overline{z}_{2}z_{3}\overline{z}_{4}\\  z_{1}{z}_{3}{z}_{4}
\end{array}\right]+t_{4}\left[\begin{array}{cc} {z}_{2}^{2}z_{3}\\ \overline{z}_{1}^{2}{z}_{4}
\end{array} \right]+ t_{5}\left[ \begin{array}{cc} \overline{z}_{1}z_{3}{z}_{4}\\
 \overline{z}_{2}\overline{z}_{3}{z}_{4} \end{array} \right]\\
 & &+t_{6}\left[
\begin{array}{cc} \overline{z}_{1}z_{2}\overline{z}_{3}\\ \overline{z}_{1}\overline{z}_{2}\overline{z}_{4}
\end{array} \right]+t_{7}\left[ \begin{array}{cc}
{z}_{1}\overline{z}_{3}^{2}\\  {z}_{2}\overline{z}_{4}^{2}
\end{array} \right]+t_{8}\left[ \begin{array}{cc} {z}_{1}\overline{z}_{2}{z}_{4}\\
 {z}_{1}{z}_{2}\overline{z}_{3} \end{array}
\right]+t_{9}\left[ \begin{array}{cc}
{z}_{1}^{2}{z}_{2}\\  \overline{z}_{1}{z}_{2}^{2}
\end{array} \right]
\end{array}
\]
Since all the maximal isotropy subgroups have positive indices, we check
how the condition of regularity on its zero set is satisfied by $Q$. Let
$z_{j}=x_{j}+iy_{j}$, $j=1,2,3,4$, and $\Sigma_{1}=\Z_4(b)$, $\Sigma_{2}=\Z_4(-b)$, and
$\Sigma_{3}=\Z_2(-b_{2})$ and letting $Q^{\Sigma_{i}}$ denote again the cubic equivariants
restricted to the respective fixed point subspaces. We have,
\[
Q^{\Sigma_{1}}(x_{1},x_{3},t)=(t_{1}+t_{2})x_{3}^{3}+(t_{3}+t_5+t_7)x_{1}x_{3}^{2}+(t_{4}+t_{6}+t_8)
x_{1}^{2}x_{3}+t_{9}x_{1}^{3}
\]
\[
Q^{\Sigma_{2}}(x_{1},x_{3},t)=(-t_{1}+t_{2})x_{3}^{3}+(t_{3}-t_5+t_7)x_{1}x_{3}^{2}+(t_{4}-t_6+t_8)
{x}_{1}^{2}x_{3}-t_{9}x_{1}^{3}
\]
and letting $y=(y_1,y_2,y_3,y_4)$ we have $(-i)Q^{\Sigma_{3}}(y,t)$
\[
\begin{array}{cl}
=&t_{1}\left[ \begin{array}{cc}  -y_{4}^{3}\\ y_{3}^{3} \end{array}\right]+t_{2}\left[ \begin{array}{cc} y_{3} y_{4}^{2}\\ -{y}_{3}^{2}y_{4} \end{array}
\right]+t_{3}\left[ \begin{array}{cc} -y_{2}y_{3}y_{4}\\ -y_{1}{y}_{3}{y}_{4} \end{array}
\right]+t_{4}\left[ \begin{array}{cc}-{y}_{2}^{2}y_{3}\\y_{1}^{2}y_{4}
\end{array}\right]+ t_{5}\left[\begin{array}{cc} y_{1}y_{3}{y}_{4}\\ -y_{2}y_{3}y_{4}\end{array} \right]\\
\\
&+t_{6}\left[\begin{array}{cc}-y_{1}y_{2}y_{3}\\ y_{1}y_{2}y_{4}\end{array}
\right]+t_{7}\left[\begin{array}{cc}-{y}_{1}y_{3}^{2}\\
 -{y}_{2}y_{4}^{2} \end{array} \right]
+t_{8}\left[\begin{array}{cc}  {y}_{1}y_{2}y_{4}\\ y_{1}y_{2}y_{3}
\end{array} \right]+t_{9}\left[\begin{array}{cc}
-y_{1}^{2}y_{2}\\  y_{1}y_{2}^{2} \end{array}\right].
\end{array}
\]
The linearization $dQ^{\Sigma}$ for the
maximal isotropy subgroups $\Sigma_1$ and $\Sigma_2$ are respectively
\[
dQ^{\Sigma_{1}}=\left[(t_{3}+t_{5}+t_{7})x_{3}^{2}+2(t_{4}+t_{6}+t_{8})x_{1}x_{3}+3t_{9}x_{1}^{2},
3(t_{1}+t_{2})x_{3}^{2}+2(t_{3}+t_{5}+t_{7})x_{1}x_{3}+(t_{4}+t_{6}+t_{8})x_{1}^{2}
\right]
\]
and
\[
dQ^{\Sigma_{2}}=\left[(t_{3}-t_{5}+t_{7})x_{3}^{2}+2(t_{4}-t_{6}+t_{8})x_{1}x_{3}-3t_{9}x_{1}^{2},
3(-t_{1}+t_{2})x_{3}^{2}+2(t_{3}-t_{5}+t_{7})x_{1}x_{3}+(t_{4}-t_{6}+t_{8})x_{1}^{2}
\right]
\]
These two cases can be analyzed as in Example~\ref{ex:D6} and one can show that
${\cal E}_G\subset \partial {\cal E}_{\Sigma_1}$ and ${\cal E}_{G}\subset \partial {\cal E}_{\Sigma_2}$. We leave the
details to the reader. We do the case of $\Sigma_3$ since the details of the computations are slightly different.
The linearization of $Q^{\Sigma_3}$ is
\[
dQ_{3}^{\Sigma_{3}}(y)=i\left[C_1(y),C_2(y),C_3(y),C_4(y)\right]
\]
where
\[
C_1(y)=\left[\begin{array}{c}
t_{5}y_{3}y_{4}-t_{6}y_{2}y_{3}-t_{7}y_{3}^{2}+t_{8}y_{2}y_{4}-2t_{9}y_{1}y_{2}\\
-t_{3}y_{3}y_{4}+2t_{4}y_{1}y_{4}+t_{6}y_{2}y_{4}+t_{8}y_{2}y_{4}+t_{9}y_{2}^{2}
\end{array}
\right]
\]
\[
C_2(y)=\left[\begin{array}{c}
-t_{3}y_{3}y_{4}
-2t_{4}y_{2}y_{3}-t_{6}y_{1}y_{3}+t_{8}y_{1}y_{4}-t_{9}y_{1}^{2}\\
-t_{5}y_{3}y_{4}
+t_{6}y_{1}y_{4}-t_{7}y_{4}^{2}+t_{8}y_{1}y_{3}+2t_{9}y_{1}y_{2}
\end{array}
\right],
\]
\[
C_3(y)=\left[\begin{array}{c}
t_{2}y_{4}^{2}-t_{3}y_{2}y_{4}-t_{4}y_{2}^{2}-t_{6}y_{1}y_{3}-2t_{7}y_{1}y_{3}-2t_{7}y_{1}y_{3}\\
3t_{1}y_{3}^{2}-2t_{2}y_{3}y_{4}
-t_{3}y_{1}y_{4}-t_{5}y_{2}y_{4}+t_{8}y_{1}y_{2}
\end{array}
\right],
\]
\[
C_4(y)=\left[\begin{array}{c}
-3t_{1}y_{4}^{2}+2t_{2}y_{3}y_{4}
-t_{3}y_{2}y_{3}+t_{5}y_{1}y_{3}+t_{8}y_{1}y_{2}\\
-t_{2}y_{3}^{2}-t_{3}y_{1}y_{4}+t_{4}y_{1}^{2}
-t_{5}y_{2}y_{3}+t_{6}y_{1}y_{2}-2t_{7}y_{2}y_{4}
\end{array}
\right].
\]
We now show that $Q^{\Sigma_3}$ is regular on its zero set for an open and
dense set of values of $(t_1,\ldots,t_9)\in \Real^{9}$.
Let $M_{ij}(y)$ be the $2\times 2$ subdeterminants obtained with columns $C_i$ and $C_j$
where $i<j$. Consider the homogeneous mapping of degree four $G:\Real^{4}\to \Real^{4}$
given by
\[
G(y)=(M_{12}(y),M_{13}(y),M_{14}(y),M_{23}(y)).
\]
From van der Waerden~\cite{vdW} Section 16.5, a set of $n$ homogeneous polynomials in $n$ variables
with complex coefficients has a nontrivial common zero if and only if a resultant
$R$ in the coefficients of the homogeneous polynomials vanishes. The resultant $R$ associated
to $G(y)$ is a function of $(t_1,\ldots,t_9)$. Therefore, $G(y)=0$ for some nonzero $y$ implies
$R(t_1,\ldots,t_9)=0$. A straightforward calculation with a symbolic algebra package shows that the
components of $G$ are distinct and linearly independent, therefore, $R(t_1,\ldots,t_9)=0$ is a proper algebraic subset of $\Real^{9}$; that is,
the complement of $R^{-1}(0)$ is an open and dense set. We choose $(t_1,\ldots,t_9)\in \Real^{9}\setminus R^{-1}(0)$ so that $Q^{\Sigma_3}$ is necessarily regular on its zero set. We choose $t\in R_G$ such that $(t_1,\ldots,t_9)\in \Real^{9}\setminus R^{-1}(0)$.
Since, $s(\Sigma_3)>0$ and $Q^{\Sigma_3}$ has odd degree, then $Q^{\Sigma_3}$ vanishes for some nonzero $y^{*}$. The remainder of the argument using Theorem~\ref{BMS} is similar to the one used in Example~\ref{ex:D6}.
Therefore, ${\cal E}_{G}\subset \partial {\cal E}_{\Sigma_3}$ and we can conclude that if a map in $C_{G}^{\infty}(V,W)$ is $G$-transverse to
$0\in W$ at $0\in V$, then the zero set in a neighborhood of the origin contains smooth $1$-dimensional submanifolds of points with isotropy
types $\Sigma_1$ and $\Sigma_2$ and smooth $2$-dimensional submanifolds of points with isotropy type $\Sigma_3$.
\end{example}

\subsection{Third Case}\label{sec:case3}
For the case $0<\delta<r$, we use a Lyapunov-Schmidt type reduction procedure.
Restrict $W$ to a subspace $W'$ which has $\delta$ irreducible representations $U$.
We decompose $W=W'\oplus W'^{\perp}$ with respect to a $G$-invariant inner product and define
\[
F^{1}:V\times\Real^{l}\rightarrow W'\quad
F^{2}:V\times\Real^{l}\rightarrow W'^{\perp}.
\]
The mapping $F^{1}$ satisfies the hypothesis of Theorem~\ref{thm:ift}
and so the implicit function theorem is applied. Let $V=U\oplus V'$ where
$U$ is $G$-isomorphic to $W'$, then there exist neighborhoods of $0$ in $U$
and $V'$ and a smooth $G$-equivariant mapping $\phi:V'\to U$ such that
$F^{1}(\phi(v'),v')\equiv 0$. Note that
\[
d_{v'}F^{1}(\phi(v'),v')|_{v'=0}=d_{u}F^{1}(0,0)d_{v'}\phi(0)+d_{v'}F^{1}(0,0)=0
\]
where $d_{v'}F^{1}(0,0)=0$ since $V'$ and $W'$ do not share any irreducible representations.
Because $d_{u}F^{1}(0,0)$ is a nonzero matrix (consequence of Schur's lemma), then
$d_{v'}\phi(0)=0$. By Theorem~\ref{thm:ift} ${\cal E}_{G}\subset \partial {\cal E}_{\Sigma}$ for
all maximal isotropy subgroups $\Sigma$ with $s(\Sigma;V,W')>0$.

We now substitute the solution of $F^{1}$ into $F^{2}$ and obtain a mapping
\[
\tilde{F}^{2}:V'\to W'^{\perp}
\]
where $d_{v'}F^{2}(0)=0$ since $d_{v'}\phi(0)=0$. Now, $V'$ and $W'^{\perp}$ do not share
irreducible representations and so we need to study $\tilde{F}^{2}$ using the method outlined
in section~\ref{sec:case2}. That is, for all maximal isotropy subgroups $\Sigma$ with $s(\Sigma;V',W'^{\perp})>0$
we need to find out if ${\cal E}_{G}\subset \partial {\cal E}_{\Sigma}$ for the mapping $\tilde{F}^{2}$.
For those $\Sigma$ which do, then we have an inclusion ${\cal E}_{G}$ in $\partial {\cal E}_{\Sigma}$ for
the full map $F$. The next example illustrates this approach.

\begin{example}\label{ex:F-group-2}
Let $V=V_3^2\times V_2$ and $W=V_2^2$. Then $W'=V_2$ and
\[
F^{1}:V\times\Real^{l}\rightarrow V_2\quad
F^{2}:V\times\Real^{l}\rightarrow V_2.
\]
\end{example}
Note that we have the same three maximal isotropy subgroups $\Sigma_i$ $i=1,2,3$ as in Example~\ref{ex:F-group} and $s(\Sigma_i)=1$ for $i=1,2$ and $s(\Sigma_3)=2$.

Using the implicit function theorem, there exists a smooth
$G$-equivariant mapping $\phi:V_3^2\to V_2$ near $0$ which solves $F^{1}=0$.
We substitute $v_2=\phi(u_3,v_3)$ in the second mapping where $(u_3,v_3)\in V_3^2$ and $v_2\in V_2$.
Since $F^{1}$ has no quadratic $G$-equivariant terms, it is easily shown by implicit differentiation that
second derivatives of $\phi$ evaluated at $0$ vanish. Therefore, as expected, the lowest degree terms in
$F^{2}(v_1,v_3,\phi(u_3,v_3))$ are cubic. Let $Q^{2}$ denote the cubic degree homogeneous truncation of $F^{2}$,
then
\[
Q^{2}(u_3,v_3,\phi(u_3,v_3))=d_{v_2}F^{2}(0,0,0)\phi(u_3,v_3)+Q_3(u_3,v_3,0).
\]
where $Q_3(u_3,v_3,v_2)$ denotes the cubic $G$-equivariant terms in $F^{2}(u_3,v_3,v_2)$.
Note that the space of cubic $G$-equivariant mappings in $C_{G}^{\infty}(V_3^2,V_2)$ is six-dimensional.
Let $u_3=(z_1,z_2)$ and $v_3=(z_3,z_4)$ the generators are:
\[
\begin{array}{c}
\left[\begin{array}{c} z_3^2 \ov{z_4}\\ z_3 z_4^2\end{array}\right],
\left[\begin{array}{c} \ov{z}_2 z_3^2 \\ z_1 z_4^2\end{array}\right],
\left[\begin{array}{c} z_1 z_3 \ov{z}_4 \\ z_2 z_3 z_4\end{array}\right],
\left[\begin{array}{c} z_1 \ov{z}_2 z_3 \\ z_1 z_2 z_4\end{array}\right],
\left[\begin{array}{c} z_1^2 \ov{z}_4 \\ z_2^2 z_3\end{array}\right],
\left[\begin{array}{c} z_1^2 \ov{z}_2 \\ z_1 z_2^2\end{array}\right].
\end{array}
\]
Calculations similar to the ones done in Example~\ref{ex:F-group} show that
${\cal E}_{G}\subset \partial {\cal E}_{\Sigma_i}$ for every isotropy subgroup
$\Sigma_i$ since we have $s(\Sigma_i;V_3^2,V_2)=1$ for $i=1,2$ and $s(\Sigma_3;V_3^2,V_2)=2$.

\section{More Questions}\label{sec:questions}
Question~\ref{Q:subgroup} in all its generality is more complicated than the problems we have investigated in this paper. The main problem is the characterization of the restriction of $F$ to a subgroup. This problem is not well understood in general since hidden symmetries~\cite{GMS-hidden} and deficiencies~\cite{Stewart-Dias2000} can modify the expected structure of the equivariant mapping on the fixed-point subspace. A positive answer to Question~\ref{Q:subgroup} leads to the following question.
\begin{question}
We know that inclusion of strata is transitive. So that if ${\cal E}_{G}\subset \partial {\cal E}_{\Sigma}$ and ${\cal E}_{\Sigma}\subset \partial {\cal E}_{\Sigma'}$, then ${\cal E}_{G}\subset \partial {\cal E}_{\Sigma'}$.
If Question~\ref{Q:subgroup} has a positive answer, a necessary condition for these inclusions is that $s(G)<s(\Sigma)<s(\Sigma')$. Is it a sufficient condition?
\end{question}
Finally, we consider the case of continuous groups
\begin{example}\label{O2ex}
The standard $\Otwo$ action on $V=W=\C$ is given by
$\kappa.z=\ov{z}$ and $\theta.z=e^{i\theta}z$. It is easy
to check (see~\cite{GSS88}) that $F(z,t)=tz$. There are two isotropy
types: $\Otwo$ and $\Z_2(\kappa)$ with
\[
{\cal E}_{\Otwo}=\{(z,t)|z=0,t\neq 0\} \quad \mbox{and}\quad
{\cal E}_{\Z_2(\kappa)}=\{(z,t)|z\neq 0, t=0\}.
\]
Here $\dim {\cal E}_{\Otwo}<\dim {\cal E}_{\Z_2(\kappa)}$ since
$s(\Otwo)=s(\Z_2(\kappa))=0$, but
${\cal E}_{\Otwo}\not\subset \partial {\cal E}_{\Z_2(\kappa)}$.
\end{example}

\begin{question}
How does one study the inclusions ${\cal E}_{G}$ into $\partial {\cal E}_{\Sigma}$ when $G$ is a continuous group. We see from Example~\ref{O2ex} that the inequality of indices is not a necessary condition for inclusion. The dimension of the normalizer certainly plays a role here and the question of group orbits of zeros also arises.
\end{question}

Future work on this topic certainly includes the questions brought up above. More importantly, we would like to understand whether Theorem~\ref{thm:ift} can be extended to the cases $\delta=0$ and $0<\delta<r$ by using a general argument along the lines of the one provided for Example~\ref{ex:D6}, Example~\ref{ex:F-group} and Example~\ref{ex:F-group-2}. One possible obstacle to the generalization is to find out whether the mapping $Q$ of lowest degree equivariant generators depends explicitly on all the variables parametrizing fixed point subspaces of maximal isotropy subgroups $\Sigma$ with positive index. Let $\{x_1,\ldots,x_n\}$ be the coordinates of an orthogonal basis of $\Fix_{V}(\Sigma)$ and suppose that $Q$ does not depend explicitly on $x_1$, then for $x_1\neq 0$, $Q(x_1,0,\ldots,0)=0$ and automatically $Q$ is not regular on its zero set for all values of $t$. Therefore, in order to use Theorem~\ref{BMS} in a general way, we need to gain an understanding of the lowest degree equivariant generators for $C_{G}^{\infty}(V,W)$ from which we can find out whether the problem of non-regularity of $Q$ on its zero set persists. The fact that Theorem~\ref{thm:faith1} restricts the problem to $V$ faithful may be crucial in pursuing this issue. There may also be other ways of generalizing Theorem~\ref{thm:ift} to the other cases of $\delta$ without using Theorem~\ref{BMS} and this would be an interesting development.

\vspace*{0.25in} \noindent {\Large\bf Acknowledgements}

\vspace*{0.2in} This research is partly supported by the Natural Sciences and Engineering
Research Council of Canada in the form of a Discovery Grant (PLB) and an Undergraduate Student Research
Award (MH).


\begin{thebibliography}{M}

\bibitem{Biers77} E.~Bierstone. General Position of Equivariant Maps.
{\em Trans.~AMS} {\bf 234(2)} (1977) 447--466.

\bibitem{BLR08} P-L.~Buono, J.S.W~Lamb, M.~Roberts. Bifurcation and branching of
equilibria in reversible equivariant vector fields. {\em Nonlinearity} {\bf 21}
(2008), 625--660.

\bibitem{BMS83} M.~Buchner, J.~Marsden and S.~Schecter. Applications
of the Blowing-Up Construction and Algebraic Geometry to Bifurcation
Problems. {\em J.Diff.Eq} {\bf 48} (1983), 404--433.

\bibitem{Field77} M.J.~Field. Transversality in $G$-manifolds.
{\em Trans.~AMS} {\bf 231(2)} (1977), 429--450.

\bibitem{field-rich1} M.J.~Field and R.~Richardson. Symmetry-Breaking and the
Maximal Isotropy Subgroup Conjecture for Reflection Groups.
{\em Arch. Rat. Mech. Anal.\/} {\bf 105} (1989), 61--94.

\bibitem{Field96}  M.J.~Field. {\em Symmetry-Breaking for Compact Lie
Groups.} Memoirs of the AMS. {\bf 574} (1996).

\bibitem{Field07} M.J. Field. {\em Dynamics and Symmetry}
ICP Advanced Texts in Mathematics vol 3, London: Imperial College Press, 2007.

\bibitem{furter} J.-E.~Furter, A.M.~Sitta and I.~Stewart. Singularity
theory and equivariant bifurcation problems with parameter symmetry.
{\em Math. Proc. Camb. Phil. Soc.} {\bf 120} (1996), 547--578

\bibitem{GWPL} C.G.~Gibson, K.~Wirthm\"uller, A.A. du Plessis and
E.J.N~Looijenga. {\em Topological Stability of Smooth Mappings}.
Lecture Notes in Mathematics v. 552, Springer-Verlag, New-York, 1976.

\bibitem{GMS-hidden} M. Golubitsky, J. E. Marsden and D. Schaeffer, Bifurcation problems
with hidden symmetries. In: {\em Partial differential equations and dynamical systems}
(ed. W. E. Fitzgibbon), Research Notes in Mathematics {\bf 101}, Pitman, San Francisco, 1984.

\bibitem{GS} M.~Golubitsky and D.G.~Schaeffer. A discussion of symmetry
and symmetry breaking. In: {\em Singularities, Part 1 (Arcata, Calif., 1981)
Proc. Sympos. Pure Math.}, {\bf 40}, Amer. Math. Soc., Providence, RI, 1983.

\bibitem{GSS88} M.~Golubitsky, I.~Stewart, D.G.~Schaeffer.
{\em Singularities and Groups in Bifurcation Theory: Vol.II}. Appl.
Math. Sci. {\bf 69}, Springer-Verlag, New York, 1988.

\bibitem{HL92} I. Hambleton and R. Lee. Perturbation of equivariant moduli spaces. {\em Math. Ann.}
{\bf 293} (1992), 17--37.

\bibitem{JL} G. James and M. Liebeck. {\em Representations and Characters of Groups}.
Cambridge University Press. Cambridge, 1993.

\bibitem{Lim-McC1} C.C.~Lim and I-H.~McComb. Time-reversible and equivariant
pitchfork bifurcation. {\em Physica D} {\bf 112} (1998), 117--121.

\bibitem{Michel72} L. Michel. Nonlinear group action: Smooth actions of compact Lie groups on manifolds.
In: {\em Statistical Mechanics and Field Theory} (R.N. Sen and C. Weil, Eds), Israel University Press, Jerusalem,
(1972), 133--150.

\bibitem{Ruelle73} D. Ruelle. Bifurcations in the presence of a symmetry group.  {\em Arch. Rational Mech. Anal.}
  {\bf 51},  (1973), 136--152.

\bibitem{Sattinger} D.H.~Sattinger. {\em Branching in the presence of
a symmetry group}. {\em CBMS-NSF Conference Notes} {\bf 40} SIAM,
Philadelphia. 1983

\bibitem{Stewart-Dias2000} I.~Stewart and A.P. Dias. Hilbert series for equivariant mappings restricted to invariant
hyperplanes. {\em J. Pure Appl. Algebra} {\bf 151}, (2000), 89--106.

\bibitem{vdW} B.L.~van der Waerden. {\em Algebra, vol 2}, Frederick Ungar
Publishing Co., New York, 1970.
\end{thebibliography}
\end{document}